\documentclass[runningheads]{llncs}
\usepackage{graphicx}
\usepackage{color}

\usepackage{diagrams}
\usepackage[all]{xy}
\usepackage{tikz-cd}
\usetikzlibrary{arrows}
\tikzset{commutative diagrams/.cd,arrow style=tikz,diagrams={>=latex'}}
\usepackage{graphicx}
\usepackage{multicol}

\usepackage{soul}
\usepackage{todonotes}

\usepackage{wrapfig}

\begin{document}
\emergencystretch 3em

\title{Shall We (Math and) Dance?
}
%
%
\author{Maria Mannone\inst{1}
\and
Luca Turchet\inst{2}
}

\authorrunning{M. Mannone and L. Turchet}
%
\institute{Department of Mathematics and Informatics, University of Palermo, Italy, \email{mariacaterina.mannone@unipa.it, manno012@umn.edu}\footnote{M. M. is an alumna of the University of Minnesota, USA.}\and
Department of Information Engineering and Computer Science, University of Trento, Italy, 
\email{luca.turchet@unitn.it}
}

\maketitle              
\begin{abstract}


Can we use mathematics, and in particular the abstract branch of category theory, to describe some basics of dance, and to highlight structural similarities between music and dance? We first summarize recent studies between mathematics and dance, and between music and categories. Then, we extend this formalism and diagrammatic thinking style to dance.

\keywords{$2$-categories \and Music \and Dance.}
\end{abstract}
\section{Introduction: Why mathematics for dance}

Joining the abstraction of mathematics with expressivity and passion in dance is possible. It means rationally exploring basic features of dance, and appreciating the flexibility of modern mathematics, for speculative investigation and practical purposes.
A joint approach between music, dance, and mathematics would involve rational thinking upon the arts, as well as the `translations' of ideas and transformational mechanisms between disciplines in a STEAM (Science, Technology, Engineering, Art, and Math) attitude.
 
`Concrete' applications can also involve pedagogy in one or more of these disciplines, giving amusing examples of applications of mathematical concepts and procedures, and investigating some hidden, theoretical roots of artistic practices. In particular, in the case of math and dance, an interplay between rationality and intuition may even help people learn dancing. This is the case of a software that helps tango learners to find the main pulse via a mathematical analysis of the underlying musical structures \cite{amiot}. A way to compare underlying rhythms for dance is suggested in E. Amiot's book  \cite{amiot_book}, where Fourier discrete transforms are used to give an idea of the distribution of durations of notes and rhythm patterns.

The world of dance includes different styles and techniques, each of them with inner symmetries. These symmetries may be thought of 
throught the lens of mathematics. An overview of mathematics for different dance styles is proposed in \cite{wasilewska}, with a focus on the geometry of figures (`poses').
The geometry of choreographies also includes connecting movements between poses of the group of dancers.
One can analyze the geometry of poses of a single dancer, but a more complete investigation should also involve the analysis of the geometry generated by the entire group of dancers on stage. The audience of a ballet, and, in general, of a dance show, pays attention differently to dancers, dancers' groups, and parts of them. To investigate the relationship between attention and single dancers 
 the concept of the {\em center of attention} has been used, an element already used in visual arts. It can be studied geometrically and statistically as a trajectory \cite{schaffer}.
In analogy with the center of mass that connects form and mass, the center of attention has been connected with the ideal `mass' of dancers, with the concept of `center of attention mass.'
The center of attention mass of an ensemble connects physics, group dancing schemes, and rhythm patterns. It is evaluated assigning weights to dancers' bodies based ``on the type of movement performed and how likely the moves are to attract the audience's attention'' \cite{wasilewska}.

A scientific perspective on choreography does not involve only analysis and understanding of dance practice: it can suggest new ideas sparkling creativity.
Furthermore, specific geometric forms can be used as bases for choreography and other artistic applications; this is the case of the truncated octahedron\footnote{An Archimedean solid.} cited in \cite{borkovitz}, and the `Apollonian circles'\footnote{Families of circles where every circle of a family intersects every circle of the other family orthogonally.} from {\em The Daughters of Hypatia} by Karl Schaffer \cite{schaffer}. 

Dance involves figures as events in time continuously connected, and time -- and expressivity -- is shaped by music. Music does also influence not only the {\em when} of figures and movements but also their {\em how}, their style. Thus, a complete study of dance cannot prescind from the role of music. In fact, music and dance are tightly connected in several cultures \cite{patel}.


One of the elements connecting music and dance is {\em gesture}, which influences expressive parameters of music. Gestures also contain information about pulse. The conductor communicates pulse and rhythm to orchestral musicians, as well as overall style and expressivity. Each musician performs specific gestures, according to the technique of the specific instrument he or she plays. Dancers do not take pulse directly from the conductor, but from performed music, via a `filtering' operation; see \cite{amiot}. Dancers' movements and gestures, especially in ballet, are linked with music and thus, indirectly, with musicians' gestures that produce that music.

We can wonder if it is possible to describe within some unitary vision the gestures of conductor, musicians, and dancers. To the best of our knowledge, there is not such a theory yet, but only a collection of experiences and case studies. A recent mathematical overview of the similarities between conducting and orchestral musicians' gestures uses the formalism of category theory \cite{mannone_jmm}. In fact, the power of abstraction of categories allows for the schematization of similarities and similar transformations between different objects, 
relaxing the condition of `equality' in favor of `up to an isomorphism.'

\begin{figure}
\centering
\includegraphics[width=12cm]{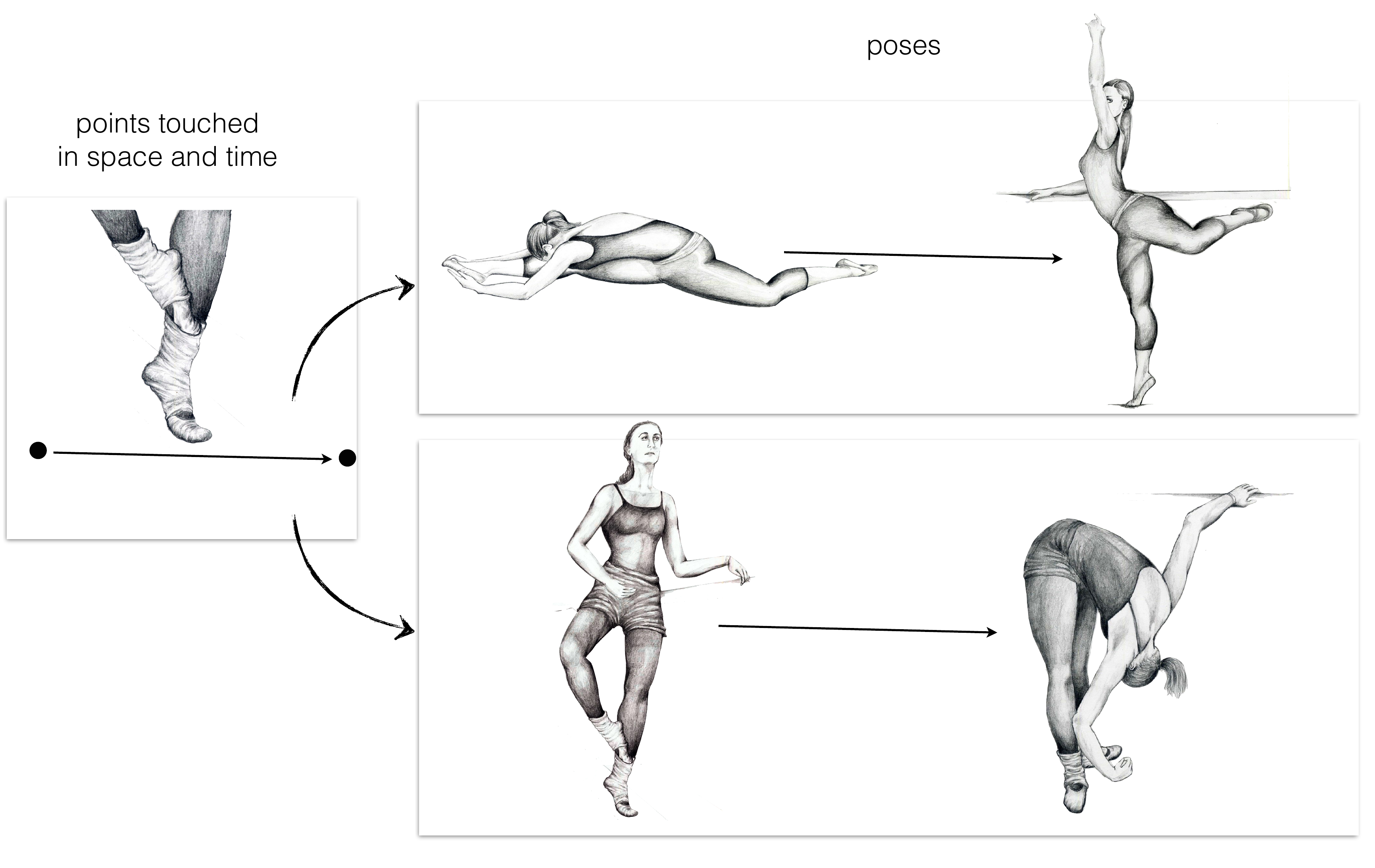}
\caption{A mapping from a simple diagram (with two points and an arrow) representing points in space and time touched by the foot of a dancer, and two different realizations of this scheme. The arrows between the two images of each pair of dancers represent the continuous movements to connect the two poses. 
The mapping consists of two (topological) functors
from points in space-time (first category) to dance figures and movements (second category). This is a category of gestures, where objects are gestures, and morphisms are hypergestures.  We are working with curves of gestures, that is, a particular instance of hypergestures. In the second category, the points are the figures (poses), and the arrows are the movements. We are using functors and not only functions because we need to map not only objects from a domain to a codomain but also the transformations between them, and functors map objects into objects, and morphisms into morphisms.
Dancers are not usually aware of this, but the mathematics behind their art is the definition of gestures on topological categories \cite{mazzola_andreatta,arias}.
Drawings by Maria Mannone.
}
\label{ballerine}
\end{figure}

The flexibility of categories and, more in general, of diagrammatic thinking is exemplified by the mathematical definition of {\em gesture}. While trying to sketch a general theory of music and dance, it is helpful to highlight common elements between different dance styles. Arrows and diagrams in categories can connect the discussed existing studies with current research in mathematical music theory. Figure \ref{ballerine} shows an 
intuitive and yet precise application of categorical formalism to dance. This schematization is inspired by the definition of musical gestures proposed in \cite{mazzola_andreatta}, that makes use of dance as a metaphor to understand music.

In this paper, we briefly summarize mathematical theory of musical gestures \cite{mazzola_andreatta,mannone_jmm,mannone_knots,arias} and basics of categories, and we show their possible application to dance. This approach may be useful to investigate formal and cognitive studies about dance \cite{patel}.

\section{Categories, Music, and Dance}\label{categories_dance}

Category theory can be applied to music \cite{mazzola_topos,fiore_noll}, included musical gestures \cite{mazzola_andreatta,arias,mazzola_topos,jed}. A category is given by objects (points or, to avoid a set-theoretic feeling, {\em vertices}) and morphisms between them (arrows), verifying associative and identity properties. We can define arrows between categories ({\em functors}) bringing objects of a category into objects of the other one, and morphisms of a category into morphisms of the other one. Considering categories themselves as vertices (to be connected via functors), we can easily build nested structures.\footnote{We can define $2$-categories and $n$-categories.}
	
In \cite{mazzola_andreatta}, a gesture is defined as a collection of vertices in space and time connected by a system of continuous curves.\footnote{More precisely, an abstract (and oriented) diagram is mapped into points and paths in a topological space.}
Musical notes (vertices) are like `the points touched by a dancer while moving continuously in the space' (arrows) \cite{mazzola_andreatta}.
This definition can be extended to visual arts, $n$-categories, and recursive gestures \cite{mannone_knots}.
Musical gestures allow one to shape physical sound parameters for expressive reasons: a purely `technical' gesture such as pressing a piano key can be shaped into, let us say, a delicate, caressing gesture.\footnote{Conversely, the violent bow's movements of strings in the movie {\em Psycho} by Hitchcock well evoke the knife hitting.}


The most basic gesture is probably {\em breathing}, that can be seen in general as a couple of (inverted) curves between two points --- a `tense' state and a `relaxed' state. In piano playing, we have the general idea of preparation and key-pressing. Breathing, in music and dance, is the starting point of both basic technique and expressive motion.\footnote{According to \cite{karin}, ``in the dance studio, conscious use of breathing patterns can enhance the phrasing and expressivity in movement.''} 
 From conductor's hands to dancers' feet, all feelings are conveyed into muscular tension/release and drawing curves in space and time.  
`Breathing' reminds of diaphragm movements for singers and wind players. More generally, breathing reminds of arsis and thesis, preparation and beat (conducting), inspiration and exhalation, dissonance and resolution in harmony. The dualism tension-relaxation is also relevant for dance \cite{charnavel}. In general, we can find {\em similarities} of gestures on different instruments --- e.g., the necessary variations of movements to play a {\em forte} on violin and on piano \cite{mannone_jmm}; these similarities also deform the basic metric gestures of a conductor.


In dance, we have poses/figures (vertices) and movements (arrows) connecting them. Movements can be seen as vertices, and transformations of these movements as arrows. Thus, both in music and dance we can categorically investigate basic technical gestures and their expressive deformations, as well as their compositions and groupings. Groupings and hierarchies of musical structures \cite{lerdhal,mazzola_topos} and gestures \cite{mannone_knots} have a correspondence for dance: basic feelings represented via expressive gestures can be connected within a whole story, letting dance acquire a narrative dimension.

\section{Diagrammatic Details}\label{diagrammatic_details}

Gestures\footnote{The category of gestures has gestures as objects and gestures of gestures (hypergestures) as the morphisms between them. The composition of hypergestures (paths) is associative up to a path of paths \cite{mannone_knots}.} in music and dance can be seen as two categories. In dance, the objects are the positions, and the morphisms are the movements to reach them.\footnote{From \cite{charnavel}: ``The physical signal of dance consists of a change in position of the dancer body in space with respect to time.'' We can see positions as vertices and changes of positions as arrows. Also: ``movement is physically created by an infinite sequence of continuous positions in space unfolding in time.'' We can represent the transition from a position to another one via an arrow representing a morphing.}
We can define a functor from music to dance, involving metrical aspects in music to be connected with metrical aspects in dance, and expressive aspects of music with expressive aspects in dance. The `vertices' would be isolated, remarkable points of the score or short sequences in music and figures in dance, and the `arrows' would be the connecting passages in music and connecting gestures in dance.


The category {\em dance} is given by positions of the human body in space and time, called {\em figures}, connected via continuous movements:
\begin{center}

$\begin{diagram}\mbox{dance figure 1} & \rTo^{movement} & \mbox{dance figure 2}\end{diagram}$. 
\end{center}

We can easily introduce 2-categories if we consider different ways to connect two figures, whose variation, represented by the double arrow, can embody changes in speed, amplitude, and expressivity; see diagram \ref{dance_cat2}.

\begin{equation}\label{dance_cat2}
\begin{tikzcd}
\mbox{dance figure 1}
 \arrow[r,bend left=30, "movement\,1"{name=U, above}] 
  \arrow[r,bend right=30, "movement\,2"{name=D, below}]  
& 
\mbox{dance figure 2}
\arrow[Rightarrow, from=U, to=D]
\end{tikzcd}
\end{equation}

The identity corresponds to the absence of movement -- keeping a static position. Horizontal composition is well-defined: the composition of two dance movements is another dance movement. Vertical composition is also well-defined: we can transform a movement into another, and the composition of two movement transformations leads to another movement. Vertical associativity is verified if we consider equivalence classes of movements connecting dance figures.\footnote{In fact, a transformation between movements requires a homotopy, and homotopy is not associative (it requires a reparametrization). However, homotopy classes are associative.} Horizontal associativity is also verified:
similarly with (musical or generic) gestures, the composition of dance hypergestures (paths) is associative up to a path of paths. Diagrams help us compare figures and movements of different styles; see diagram \ref{dance_functor}. 
In fact, we can define, within the category `dance,' a subcategory for dance-style 1, and another subcategory for dance-style 2. For sake of simplicity, we can just define a category for style 1 and a category for style 2.
\begin{equation}\label{dance_functor}
\begin{footnotesize}
\begin{diagram}
\mbox{dance figure 1 (style 1)} & \rTo^{movement\,(style\,1)} & \mbox{dance figure 2 (style 1)}
\\ \dMapsto^{figure\,mapping} & \dMapsto^{gesture\,mapping} & \dMapsto^{figure\,mapping} 
\\ \mbox{dance figure 1 (style 2)} & \rTo^{movement\,(style\,2)} & \mbox{dance figure 2 (style 2)}
\end{diagram}
\end{footnotesize}
\end{equation}
\begin{footnotesize}
\begin{equation}\label{dance_functor_2} 
\begin{tikzcd}[swap,bend angle=45]
\mbox{dance figure 1 (style 1)}
\arrow[mapsto]{dddd}{figure\,mapping} 
    \arrow[r,bend right=20, "movement\,2\,(style\,1)"{name=D}]  
    \arrow[r,bend left=20, "movement\,1\,(style\,1)"{name=U, above}] 
&
  \arrow[Rightarrow, from=U, to=D, " "]  
\mbox{dance figure 2 (style 1)}
\arrow[mapsto]{dddd}{figure\,mapping} 
   \\  
   \\ 
   \\ & & & & 
   \\ 
   \mbox{dance figure 1 (style 2)}
  \arrow[r,bend right=20, "movement\,2\,(style\,2)"{name=Z}]
   \arrow[r,bend left=20, "movement\,1\,(style\,2)"{name=V, above}]  
&
\mbox{dance figure 2 (style 2)} 
   \arrow[Rightarrow, from=V, to=Z, " "] 
\end{tikzcd}
\end{equation}
\end{footnotesize}
We can group movements in two, three, and so on, according to tempo. Dance teachers can ask their students to start clapping hands, letting them think of musical rhythm. Then, students will start making movements with their bodies according to the given rhythm.
Mathematically, this is a mapping from the category of pulses to the category of basic dance movements; thus, we have a functor. In the category of pulses, objects are accents (beats), and arrows are time intervals between them. Beats are mapped into positions of the body in space and time, and time intervals are mapped into body movements to reach these positions; see diagram \ref{pulse_dance_2},\footnote{In the case of total improvisation (music and dance), the diagram would be commutative.} enriched with $2$-categories. Diagram \ref{dance_functor} is enriched in terms of $2$-categories in diagram \ref{dance_functor_2}. In diagram \ref{dance_functor_2}, additional arrows, not shown for reasons of graphic clarity, map $movement\,1\,(style\,1)$ into $movement\,1\,(style\,2)$, $movement\,2\,(style\,1)$ into $movement\,2\,(style\,2)$, and double arrows into double arrows. The same applies to diagram \ref{pulse_dance_2}, where $2$-cells of $2$-category pulse may represent the time variation between two consecutive pulses. We are using the formalism of $2$-functors, that are morphisms between $2$-cells\footnote{As defined in \cite{maclane}, a $2$-functor between two categories $A$ and $B$ is a triple of functions that map objects, arrows, and $2$-cells of $A$ into objects, arrows, and $2$-cells of $B$ respectively, preserving ``all the categorical structures.''} \cite{maclane}. In fact,  $2$-cells give a metaphorical idea of the different ways to connect figures within the same style, to move hands in conducting between two pulses (diagram \ref{music-to-dance_2}) for technical and expressive reasons, and so on.

\begin{footnotesize}
\begin{equation}\label{pulse_dance_2}
\begin{tikzcd}[swap,bend angle=45]
\mbox{pulse 1}
\arrow[mapsto]{dd}{} 
    \arrow[r,bend right=20, "time\,interval\,2"{name=D}]  
    \arrow[r,bend left=20, "time\,interval\,1"{name=U, above}] 
&
  \arrow[Rightarrow, from=U, to=D, " "] 
\mbox{pulse 2}
\arrow[mapsto]{dd}{} 
   \\ & & & &  
   \\ 
   \mbox{dance figure 1}
  \arrow[r,bend right=20, "movement\,2"{name=Z}]
   \arrow[r,bend left=20, "movement\,1"{name=V, above}]  
&
\mbox{dance figure 2} 
   \arrow[Rightarrow, from=V, to=Z, " "]
\end{tikzcd}
\end{equation}
\end{footnotesize}
Figures in dance can be different for women and men, for single or couple dance, for groups of three or more people, and so on. Thus, there will be different mappings from pulse to dance figures, but dance figures of all dancers, and their connecting movements, have to be coordinated. Temporal consistency is guaranteed by the action of the functor from time pulse to dance. The specific connection between dance figures and movements depends on the style of dance and the role of each dancer. The personal contribution of individual dancers is mathematically represented by additional arrows that modify dance figures and their connecting movements. Transformations between poses (i.e., movements) should take into account the expression of physical constraints: for example, the distance between dancers' bodies.
Tempo transformations play an important role: for example, a basic ternary movement can be transformed into a `valzer' movement by prolonging the first pulse.

The upper side of diagram \ref{music-to-dance_2} represents the mapping from rhythmic pulse to conducting gestures (r. h. indicates right hand),  where we consider the $2$-cell given by two different right-hand's movements connecting two right-hand's positions; the lower side, the mapping from the same rhythmic pulse to basic dance movements (with a similar remark on $2$-categories). Movements in conducting and dance are related because dancers hear the music, and they find conductor's beats through the music. If the dancers are dancing on their own without any music, they have to think tempo pulse in their own, too;\footnote{Dancers can move without any (external) music if they are able to think of their pulse and to communicate with each other via touch and non-verbal indications (and the leader should be particularly clear in such indications), especially for couple dancing; but this is an unstable situation. Permutations of roles should be clearly signaled via, for example, a variation in touch or visual communication.} if there is music, dancers can rely on its pulse. Ideally, dancers, musicians, and of course the conductor, are thinking the same tempo pulse. For the sake of completeness, we should cite the Cunningham school, where dance is not connected with the pulse of music. In this case, 
we cannot define a functor from pulse to dance.
We might describe diagram \ref{music-to-dance_2} as a simple functor from pulse to dance; however, it stresses the `generator' role of pulse for both dance and conducting.
We can imagine a `pulse' that is first an abstract thought and then is embodied in conducting gestures and dance movements. In the practical reality, dancers extract, as in a filter operation, the information about tempo pulse from music. Thus: pulse \emph{is the basis of} (is mapped into) conducting gestures, that \emph{are the basis of} orchestral playing, that \emph{is the basis of} listening to music, \emph{that helps} dancers catch pulse from music, \emph{allowing them} to move accordingly; see diagram \ref{conductor_musician_dancer}. 



\begin{footnotesize}
\begin{equation}\label{music-to-dance_2}
\begin{tikzcd}[swap,bend angle=30] 
\mbox{conductor's r.h. position 1}
    \arrow[r,bend right=20, "hand's\,movement\,2"{name=D}]  
    \arrow[r,bend left=20, "hand's\,movement\,1"{name=U, above}] 
&
  \arrow[Rightarrow, from=U, to=D, " "] 
\mbox{conductor's r.h. position 2}
   \\
   \\
   \\ & & & & 
   \\ 
   \arrow[mapsto]{uuuu}{pulse-cond.} 
\mbox{pulse 1}
\arrow[mapsto]{dddd}{pulse-dance} 
    \arrow[r,bend right=20, "time\,interval\,2"{name=D}]  
    \arrow[r,bend left=20, "time\,interval\,1"{name=U, above}] 
&
  \arrow[Rightarrow, from=U, to=D, " "] 
   \arrow[mapsto]{uuuu}{pulse-cond.}
\mbox{pulse 2}
\arrow[mapsto]{dddd}{pulse-dance} 
   \\
   \\
   \\ & & & & 
   \\ 
   \mbox{dance figure 1}
  \arrow[r,bend right=20, "body's\,movement\,2"{name=Z}]
   \arrow[r,bend left=20, "body's\,movement\,1"{name=V, above}]  
&
\mbox{dance figure 2} 
   \arrow[Rightarrow, from=V, to=Z, " "]
\end{tikzcd}
\end{equation}
\end{footnotesize}
The arrow from  conductor to dancer is dashed because the conductor does not directly give tempo to the dancer: he or she catches it from music, thus, from musicians' activity.
\begin{equation}\label{conductor_musician_dancer}
\begin{footnotesize}
\begin{diagram}
 & &\mbox{conductor}
\\  & \ldTo^{gives\,tempo} & & \rdDashto^{gives\,tempo} 
\\ \mbox{musician}  & & \rTo^{suggests\,tempo} & & \mbox{dancer}
\end{diagram}
\end{footnotesize}
\end{equation}

The mapping from tempo-pulse to dance can be formally described via a {\em pulse-dance functor}; see diagram \ref{pulse_dance_2}.
\begin{equation}\label{conductor_musician_dancer_listener}
\begin{footnotesize}
\begin{diagram}
 & &\mbox{conductor}
\\  & \ldTo^{gives\,tempo} &   & \rdDashto^{gives\,tempo} 
\\ \mbox{musician}  & & \lImplies & & \mbox{dancer}
\\ & \rdTo^{\,reaches\,via\,sound} & & \ldTo^{reaches\,via\,image\,}
\\ & & \mbox{audience}
\end{diagram}
\end{footnotesize}
\end{equation}

Overall musical character and expressivity shape expressive dance movements. For example, a sudden orchestral {\em forte} should imply a corresponding choreography variation for ballet dancers.\footnote{For tango music, this is more complicated than a simple `pulse extraction.' We can instead think of categories enriched with maps from a beat to another beat containing inner maps, that is, inner pulses. 
In salsa, there often are recognizable patterns with clave rhythms between strong beats. The pulse can depend on the specific music style, and it can be the object of future research in itself.}
We can navigate through complexity of dance, by considering comparisons between gestures of orchestral musicians within time, and gestures of dancers. Such a complex mapping can be formalized as a {\em Choreography functor}, corresponding to the choices of a choreographer that instructs ballet dancers on how to move on stage while listening to the orchestra. In folk music, this functor connects traditional instruments' playing with folk dancers' movements.
Intuitively, a choreography functor should map marked points of the musical score corresponding to groups of musical gestures into dance figures of a group of dancers, and morphisms between these points of the score (and the associate connecting groups of gestures) into `connecting' group dance movements. In this framework, natural transformations between two choreography functors would describe differences between two different choreographies based on the same music.

\begin{figure}
\centering
\includegraphics[width=6.5cm]{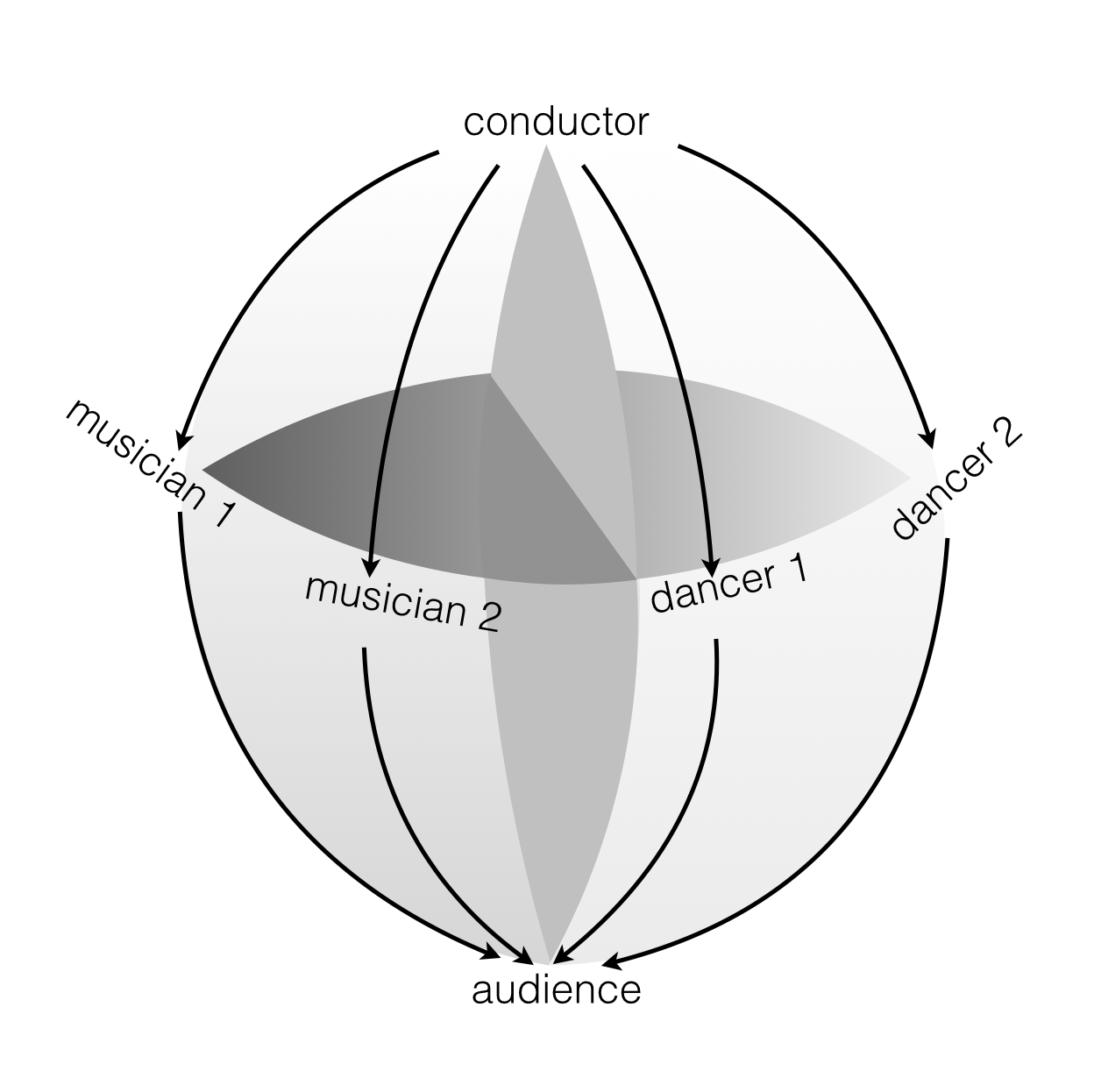}
\caption{Ellipsoid where conductor and audience can be thought of as limit and colimit, respectively.}
\label{ellipsoid}
\end{figure}

We can schematize the action of conductor, orchestral musicians, and dancers, connecting them with listeners in the audience. In diagram \ref{conductor_musician_dancer_listener}, both music played by musicians and choreography made by dancers reach the final target, the audience. Here, the conductor metaphorically plays the role of a limit, while the audience plays that of a colimit.\footnote{Limits and colimits are generalizations of products and coproducts, respectively; they are obtained the ones from the others via reversing arrows \cite{lawvere}.} Arrows appear thus inverted with respect to \cite{mannone_jmm}, but the two representations are both valid. See details in \cite{mannone_jmm} about conductor and listener. As noted by one of the reviewers,  diagram \ref{conductor_musician_dancer_listener} does not commute, but it does up to a 2-cell as indicated by the double arrow. The reason is that impressions of the public with respect to musicians and with respect to dancers can be similar, but they will never be equal. The commutativity is verified if we consider a single listener, or if we consider the equivalence class of all impressions of the public as a whole.
More pictorially, we can imagine diagram \ref{conductor_musician_dancer_listener} as an ellipsoid (Figure \ref{ellipsoid}), with the larger horizontal section representing the connections between musicians' and dancers' movements, and the vertical one, the separation between music and dance worlds.\footnote{An ellipse is used for the center of attention in \cite{schaffer}. The ellipsoid of Fig. \ref{ellipsoid} is a diagram connecting different elements. If the listener/audience recovers the pulse contained in conducting gestures, the two extremities of the ellipsoid can be joined, transforming it into a torus with a section collapsed into a point. 
}
The conductor's role appears as being more relevant for a ballet than for other dance styles. Diagrams \ref{dance_functor} and \ref{conductor_musician_dancer_listener} are given for completeness, and to connect our study with former studies in mathematical theory of musical gestures of orchestral musicians \cite{mannone_jmm}.


\section{Discussion}\label{discussion}


In this paper, we introduced the use of category theory, an abstract branch of mathematics, to the study of dance  in relationship with music. The nested structures of categories help connecting different perspectives on dance, such as the relationship with musical pulse \cite{amiot} and the geometry in dance \cite{borkovitz}.
Also, the `trajectory,' that is, the change over space/time of the center of attention as a trajectory \cite{schaffer} can be schematized via an arrow, and the comparisons between centers of attention for different choreographies, via arrows between arrows.
This would also allow for a comparison between the geometry of different choreographies, combining into a unitary vision the collection of examples proposed in literature \cite{borkovitz,schaffer}.
The geometry of choreographies can be investigated via connected categories.\footnote{A connected category $J$ is a category where, for each couple of objects $j, k\in J$, there is a finite sequence of objects $j_0,j_1,...,j_{2n}$ connecting them, that is:  $j=j_0\rightarrow j_1  \leftarrow j_2 \rightarrow... \rightarrow j_{2n-1} \leftarrow j_{2n}=k$, where both directions are allowed \cite{maclane}. Here, the morphisms are the arrows $f_i:j_i\rightarrow j_{i+1}$ or $f_i:j_{i+1}\rightarrow j_{i}$.} In the case of the tango ball, we can also think of connected musical sequences, the so-called {\em tandas}, with their different styles.

Also, one can establish a connection between dancers' positions and movements and their center of attention mass \cite{wasilewska}, and one can compare variations of this center within a choreography with variations of this center within another choreography, quantitatively characterizing choreographers' styles and different editions/productions of the same ballet. More generally, a comparison between variations of the center of attention mass would well characterize differences and analogies of choreographies for different dance styles.

We can wonder if nested structures of categories may help connecting even more layers of understanding and aesthetic appreciation for dance. We can start from the symmetry of the single, idealized human body to reach the most complex choreography. The possible steps of formal analysis would be:
\begin{enumerate}
\item the beauty\footnote{The definition of {\em beauty} is beyond the scope of this paper, and it would start a philosophical debate. We can just say that a mixture of symmetry, balance, proportion, and smoothness of movements can be overall thought of and mathematically investigated as `beauty' in dance.} of the human body as part of nature;
\item the beauty of different poses of the human body in dance, as photographical shoots;
\item the beauty of dancing movements of the human body, that is, of the connections and transitions between static poses of dance;
\item the beauty of several people dancing together: `static poses' as well as their connecting movements.
\end{enumerate}

If, on the one hand, we tend towards abstraction and generalization, on the other hand we can also wonder about the applicability of our model in concrete setups and in pedagogical frameworks, suggesting categorical developments of recent research and maybe new software applications \cite{amiot}. Both points and arrows are part of dance learning: dance students learn the poses and the way to smoothly, and expressively, reach them. The movements of each dancer have to be well-coordinated with movements of the other dancer of the couple, and with all the other dancers involved, if any. Thus, the formalism of monoidal categories\footnote{In a nutshell, a monoidal category, also called a tensor category, is a category $C$ having a bifunctor $\otimes:C\times C\rightarrow C$, that verifies pentagonal and triangular identities \cite{lawvere}. See \cite{popoff2,jed,mannone_knots} for examples of monoidal categories in music.} can catch the `simultaneity' of transformations. We need monoidal categories and not just categories because, having more dancers dancing simultaneously, with monoidal categories we can easily model their different movements with different transformations.

In fact, next developments of this research can concern an abstract description of group dancing, from couples to larger groups.
Diagram \ref{couple_dance_1} refers to monoidal categories to represent the movements of a couple of dancers. Let $A,\,B$ be two dancers, and $P^A_i,\,P^B_i$ their $i$-th poses (or `figures'), respectively. If dancer $A$ changes figure while dancer $B$ remains steady, the overall movement can be schematized by the action of the morphism $t\otimes 1$, where $t$ indicates the movement, and $1$ is the identity. Conversely, if $A$ is steady and $B$ moves, the movement will be represented by $1\otimes t'$. In general, both dancers are moving, thus we have $t\otimes t'$ in the simplified case where  $A$ is moving via $t$ and $B$ via $t'$.

\begin{equation}\label{couple_dance_1}
\begin{diagram}
P_1^A\otimes P_1^B & \rTo^{t\otimes 1} & P_2^A\otimes P_1^B
\\ \dTo^{1\otimes t'} & & \dTo^{1\otimes t'}
\\ P_1^A\otimes P_2^B & \rTo^{t\otimes 1} & P_2^A\otimes P_2^B
\end{diagram}
\end{equation}

Diagram \ref{couple_dance_2} shows the action of a `smoothness' operator $s$, that changes the movement of the first dancer --- let us say, a beginner -- from a {\em not smooth} movement to a {\em smooth} one. The complete operator, indicated in the diagram by the double arrow, is $s\otimes 1$, because it acts as an identity with respect to the second dancer, that is already moving smoothly. In diagram \ref{couple_dance_2}, the final poses are the same but the way to reach them are different. If dancers stop moving but music continues, we can describe dancers' movements via identities leading to the same positions. Also, if the leader role permutes during the show, we can use braids in monoidal categories. Finally, the formalism we developed for dance can be extended to other movements that `have to be made in a specific way': a reviewer thought of a worker building a wall, or a tourist walking while visiting a church. This might raise a question: when is a movement `dance'? We can think of formal characteristics, such as smoothness of movements and music-pulse following, but we can also think of aesthetics (and motivation) of movement in itself. According to \cite{charnavel}, ``Everyday movements are usually goal-oriented, which makes them fundamentally different from dance.''

\begin{equation}\label{couple_dance_2}
\begin{tikzcd}
P_1^A\otimes P_1^B
 \arrow[r,bend left=40, "not\,smooth\otimes smooth"{name=U, above}] 
  \arrow[r,bend right=40, "smooth\otimes smooth"{name=D, below}]  
& 
P_2^A\otimes P_2^B
\arrow[Rightarrow, from=U, to=D,"s\otimes 1"]
\end{tikzcd}
\end{equation}

\section{Conclusion}
Summarizing, we applied to dance categorical formalism formerly used for music \cite{mannone_knots,mannone_jmm}, comparing gestures in music and dance, and highlighting some of their connections.
We started from conceptual considerations in Section \ref{categories_dance}, we discussed metaphorical yet formal applications of functors and commutative diagrams to music and dance in Section \ref{diagrammatic_details}, and we finally made some references to monoidal categories in Section \ref{discussion}.
Previous work concerned mathematical description of gestures \cite{mazzola_andreatta}, connections of categorical and topological formalism with function spaces to investigate spaces of gestures \cite{arias}, and application of $2$-categories and $n$-categories to orchestral playing and to describe recursive musical gestures \cite{mannone_knots,mannone_jmm}. Former works about a formal and often mathematical description of dance topics  \cite{borkovitz,charnavel,schaffer,wasilewska} motivated us to include dance within the aforementioned diagrammatic and categorical formalism, trying to envisage general characteristics in different dance styles, and their general connection with music and musical gestures.
This paper aimed to contribute towards a unified vision of  gestures of conductor, musicians, and dancers. The mathematical reason behind that is to stress the importance of the unifying power of diagrammatic thinking to navigate within the complexity of the arts and to connect in a simple way things whose nature appears at first as deeply different. The artistic reason is to stress the importance of mutual knowledge between artistic roles, and the relevance of a well-working communication and continuous exchange via sounds, images, and gestures. We can try to approach expressivity in the arts in a simple and rational way. We can do the same thing while admiring beauty of nature via scientific investigation. This is an attempt to find keys and hidden rules of aesthetics\footnote{Categorical approaches on aesthetics of the arts \cite{kubota} and of mathematics and the arts \cite{mannone_aesthetics} can be joined with studies of nature.} and beauty, their hidden skeleta constituting the structure of a magnificent building.
Future research may include a more detailed description of particular musical genres, such as tango, and its connection with tango music \cite{amiot}. We can investigate if advanced topics in category theory do have a meaningful application in dance, and, vice versa, how essential topics in dance can be categorically described.  
Also, generative theories used to modeling natural language and dance \cite{charnavel,patel} and natural forms (with their growth processes) \cite{lindenmayer} can be compared within a categorical framework, looking for formal analogies and differences, hoping to find hidden connections between beauty in dance (and music) and beautiful natural forms, looking for some natural `roots' of beauty in the arts.
And yes, let's (math and) dance!

\section{Acknowledgments}

The authors are grateful to the mathematician, musician, and tango dancer Emmanuel Amiot for his helpful suggestions.

%
%
%


\end{document}